\documentclass[oneside,12pt]{article}

\usepackage{geometry}
\usepackage{amsmath,amssymb}
\usepackage{csquotes}
\usepackage[colorlinks=true]{hyperref}
\usepackage{authblk}
\usepackage{textcomp}
\usepackage{thmtools}
\usepackage{paralist}
\usepackage{graphicx}
\usepackage{caption}
\usepackage[labelformat=simple]{subcaption}

\usepackage[nodisplayskipstretch]{setspace}
\usepackage{mathabx}
\usepackage{mathdots}
\usepackage{booktabs}
\usepackage{multicol}
\usepackage{enumitem}
\usepackage{mathtools}

\newtheorem{theorem}{\textsc{Theorem}}[section]
\newtheorem{lemma}[theorem]{\textsc{Lemma}}

\newtheorem{remark}[theorem]{\textsc{Remark}}

\setlist[enumerate,1]{label = \arabic*),ref = \arabic*)}
\setlist[enumerate,2]{label = \emph{\alph*}),ref = \theenumi.\emph{\alph*}}

\mathchardef\ordinarycolon\mathcode`\:
\mathcode`\:=\string"8000\def\R{\mathbb{R}}
\begingroup \catcode`\:=\active
  \gdef:{\mathrel{\mathop\ordinarycolon}}
\endgroup

\renewcommand{\R}{\mathbb{R}}
\newcommand{\C}{\mathbb{C}}

\usepackage{mathtools}
\usepackage{mathrsfs}

\newcommand{\Bc}{\mathcal{B}}

\newcommand{\Rc}{\mathcal{R}}

\newcommand{\norm}[1]{\left\lVert#1\right\rVert}      
\newcommand{\abs}[1]{\left|#1\right|}                 
\newcommand{\paren}[1]{\left(#1\right)}               
\newcommand{\sparen}[1]{\left\{#1\right\}}		      
\newcommand{\h}[1]{\left\langle#1\right\rangle}       

\renewcommand{\d}{\,\mathrm{d}}						  

\DeclareMathOperator*{\argmin}{\mathrm{arg\:min}}

\renewcommand{\epsilon}{\varepsilon}
\renewcommand{\rho}{\varrho}
\renewcommand{\phi}{\varphi}

\usepackage{xcolor}
\definecolor{darkblue}{rgb}{0.1,0.1,0.65}
\definecolor{darkgreen}{rgb}{0.1,0.65,0.1}
\definecolor{darkred}{rgb}{0.65,0.1,0.1}

\definecolor{myblue}{rgb}{.95, .95, 1}
\definecolor{myred}{rgb}{0.975, .95, .95}

\title{Combining reconstruction and edge detection in computed tomography}

\author{J\"urgen~Frikel$^1$, Simon~G\"oppel$^2$, Markus~Haltmeier$^2$}

\affil{ 
\textsuperscript{1}Department of Computer Science and Mathematics, OTH Regensburg
\textsuperscript{2}Department of Mathematics,
University of Innsbruck\authorcr
 {\tt juergen.frikel@oth-regensburg.de}
 }

\date{}

\begin{document}

\maketitle

\begin{abstract}
We present two methods that combine image reconstruction and edge detection in computed tomography (CT) scans. Our first method is as an extension of the prominent filtered backprojection algorithm. In our second method we employ $\ell^{1}$-regularization for stable calculation of the gradient. As opposed to the first method, we show that this approach is able to compensate for undersampled CT data.
\end{abstract}

\section{Introduction}

Detection of edges in computed tomography (CT) scans is a challenging task because the underlying image reconstruction problem is ill-posed and, hence, even small errors in the x-ray measurements may lead to huge reconstruction errors that have to be compensated by the edge detection algorithms, cf. \cite{Natterer86}. This task becomes even more challenging, when the CT scans are generated from a small number of x-ray measurements, as it is the case, e.g., in digital breast tomosynthesis or dental CT, where the x-ray measurements can be taken only from a small number of views in a limited angular range. 
In such situations, classical reconstruction algorithms may generate characteristic reconstruction artifacts (in addition to noise amplification) and, thus, substantially degrade the image quality, cf. \cite{FrikelQuinto2013}. Performing edge detection after image reconstruction, can therefore lead to unreliable edge maps. In this article we present two methods, that allow for a stable reconstruction of edges directly from CT data and, thus, stabilizes edge detection in CT images. In this work, we focus on stable reconstruction of the  gradient of a CT scan, because it is the main ingredient in most prominent edge detection algorithms, such as the Canny algorithm, cf. \cite{CannyOrig}. 

In what follows, we model CT scans as functions $f:\R^{2}\to\R$ and the corresponding CT data by the Radon transform of $f$, which is defined as (cf. \cite{Natterer86})
\begin{equation}
\label{eq:radon transform}
	\Rc f(\phi,s) := \int_{\R} f(s\theta(\phi)+t\theta^\perp(\phi)) \d t,
\end{equation}
where $\theta(\phi)=(\cos(\phi),\sin(\phi))^{\top}$ and $\theta^{\perp}(\phi)=(-\sin(\phi),\cos(\phi))^{\top}$. Here, the value $\Rc f(\phi,s)$ represents one x-ray measurement along the x-ray path that is given by the line $L(\phi,s) = \sparen{s\theta(\phi)+t\theta^\perp(\phi): t \in \R} = \sparen{x\in\R^{2}: \h{x,\theta(\phi)}=s}$, where $\phi\in[0,\pi)$ and $s\in\R$. 

We will present two methods for stable reconstruction of the gradient of the smooth function \(f_{\epsilon}:=f\ast g_{\epsilon}\),

where 
\begin{equation}
\label{eq:gaussian}
	g_{\epsilon}(x) := \frac{1}{2\pi\epsilon^{2}}\exp\paren{-\frac{\norm{x}^{2}}{2\epsilon^{2}}}, \quad \epsilon>0.
\end{equation}

First, we derive a method that is of the same type as the famous filtered backprojection algorithm (FBP). This method  follows the spirit of  \cite{Lo2008} and yields good results whenever the CT data is well sampled \cite{Natterer01}. In our second approach, we employ ideas from compressed sensing and propose to use sparse regularization for the reconstruction of the image gradient. We show that this method leads to a more robust edge detection, especially when the data is not sampled properly.

\section{Materials and methods}

\subsection{Method 1: An FBP-type approach for calculating the gradient}
\label{sec:method1}
  
In order to calculate the smoothed partial derivatives derivatives
\begin{equation}
\label{eq:approximate derivative}
	\frac{\partial f_{\epsilon}}{\partial x_{j}} =  \frac{\partial}{\partial x_{j}}(f\ast g_{\epsilon}) = f\ast \frac{\partial g_{\epsilon}}{\partial x_{j}}, 
\end{equation}
directly form CT data, we use the well-known relations between the Radon transform, convolutions and derivatives (cf. \cite{Natterer86}) and obtain the following result.

\begin{lemma}
\label{lem:data derivative}
	Let the Gaussian $g_{\epsilon}$ be defined by \eqref{eq:gaussian}. Then, we have for $j\in\sparen{1,2}$:
	\begin{equation}
	\label{eq:radon preprocessed data}
		\Rc\left[ f\ast \frac{\partial g_{\epsilon}}{\partial x_{j}}\right] (\phi,s) =  \left[ \Rc f \ast_{s} G_{\epsilon,j}\right] (\phi,s),
	\end{equation}
	where $\ast_{s}$ denotes the convolution with respect to the second variable $s$, and
	\begin{equation}
	\label{eq:radon derivative filter}
		G_{\epsilon,j}(\phi,s) := - \frac{\theta_{j}(\phi)}{\epsilon^{3}\sqrt{2\pi}} \cdot s \cdot \exp\paren{-\frac{s^{2}}{2\epsilon^{2}}}
	\end{equation}
	with $(\phi,s)\in[0,\pi)\times\R$ and $\theta_{1}(\phi) :=\cos(\phi),~ \theta_{2}(\phi):=\sin(\phi)$.
\end{lemma}

Lemma \ref{lem:data derivative} shows that the partial derivatives \eqref{eq:approximate derivative} can be obtained by applying the inverse Radon transform $\Rc^{-1}$ to the preprocessed data \eqref{eq:radon preprocessed data}, i.e.,
\begin{equation}
\label{eq:inverse radon derivative}
	\frac{\partial f_{\epsilon}}{\partial x_{j}} = \Rc^{-1}\left[ \Rc f \ast_{s} G_{\epsilon,j} \right].
\end{equation}
Since the FBP algorithm is a regularized implementation of $\Rc^{-1}$ (cf. \cite{Natterer86}), a standard toolbox implementation could be used in practice to obtain  $\nabla f_{\epsilon}$.

Since the inverse Radon transform is given by \(\Rc^{-1} = \Bc \circ P\), where $P$ is a filtering (convolution) operator and $\Bc$ is the so-called backprojection operator (cf. \cite{Natterer86}), the filtering by $P$ can be combined with the filtering in \eqref{eq:inverse radon derivative}. In this way, we obtain a reconstruction formula of FBP-type for the derivatives \eqref{eq:inverse radon derivative}, that is stated in the following theorem.

\begin{theorem}
\label{thm:filter derivative red}
	Let $W_{\epsilon,j}$ be the function of two variables $(\phi,s)$ which is defined in the Fourier domain by
\begin{equation}
\label{eq:filter derivative rec}
	\widehat W_{\epsilon,j} (\phi,\omega) := \frac{1}{4\pi}\cdot\theta_{j}(\phi)\cdot i\cdot \omega\cdot \abs{\omega}\cdot \exp\left(\frac{-\epsilon^{2}\omega^{2}}{2}\right),
\end{equation}
where $\widehat W_{\epsilon,j} $ denotes the 1D-Fourier transform of $W_{\epsilon,j}$ with respect to the $s$-variable and $i\in\C$ is the imaginary unit. Then,

\begin{equation}
\label{eq:fbp rec grad}
	\frac{\partial f_{\epsilon}}{\partial x_{j}} = \Bc( \Rc f \ast_{s} W_{\epsilon,j}),
\end{equation}
	where $\Bc$ is the backprojection operator for the Radon transform (cf. \cite{Natterer86}).
\end{theorem}

The proof of Theorem \ref{thm:filter derivative red} follows from \eqref{eq:radon preprocessed data} together with the convolution theorem for the Fourier transform and Theorem II.2.1 in \cite{Natterer86}. It shows that the derivatives \eqref{eq:approximate derivative} can be calculated using a FBP-algorithm with angle-dependent filters that are given by \eqref{eq:filter derivative rec}.

\begin{remark}
	The accuracy of the presented method (given by \eqref{eq:inverse radon derivative} and \eqref{eq:fbp rec grad}) will strongly depend on the sampling of the Radon transform, cf. \cite{Natterer86}. If the CT data does not satisfy the sampling requirements, e.g., when the angles are sampled rather sparsely, the algorithm will produce artifacts which can substantially degrade the performance of edge detection. 
\end{remark}

\subsection{Method 2: Gradient calculation using $\ell^{1}$-regularization}
\label{sec:method2}
In order to account for the negative effects of (possible) undersampling of CT data, we propose to replace $\Rc^{-1}$ in \eqref{eq:inverse radon derivative} by a regularization method for $\Rc^{-1}$ that is able to deal with undersampled data. From the theory of compressed sensing it is well known that sparsity can help to overcome the classical Shannon sampling paradigm, \cite{CScrt}. As we are interested in recovering gradients of images, which only have large values around edges, we aim at enforcing sparsity of image gradients. Thus, we calculate the derivatives \eqref{eq:approximate derivative} approximately  via (cf. \cite{CScrt}):
\begin{equation}
\label{eq:ell1 regularization}
	\mathbf{f}_{\lambda,\epsilon}^{(j)} = \argmin_{\mathbf{f}}\norm{\mathbf{R}\mathbf{f}-\mathbf{y}\ast G_{\epsilon,j}}_{2}^{2} + \lambda\norm{\mathbf{f}}_{1},
\end{equation}
where $\lambda>0$ is regularization parameter. Note that in \eqref{eq:ell1 regularization} the bold face symbols denote the discretized versions of the corresponding continuous objects.

\section{Results}
\label{sec:results}
We implemented the methods 1 and 2 in Matlab and tested them on the x-ray CT data of a lotus root, cf. \cite{FIPSlotus}. In our experiments we converted the fan-beam data to a parallel-beam data using Matlab function \texttt{fan2para} and downsampled this data in order to simulate angular undersampling. Thereby, we used 738 equispaced samples in the $s$-variable and $36$ evenly distributed angles in $[0,\pi)$. For our implementations, we used the Matlab functions \texttt{radon} and \texttt{iradon} as numerical realizations of the Radon transform, the inverse Radon transform and the backprojection operator. For the minimization of the $\ell^{1}$-functional \eqref{eq:ell1 regularization}, we implemented the iterative soft thresholding algorithm, cf. \cite{DDD04}.

\begin{figure}[h]
     \centering
     \begin{subfigure}[t]{0.3\textwidth}
         \centering
         \includegraphics[width=\textwidth]{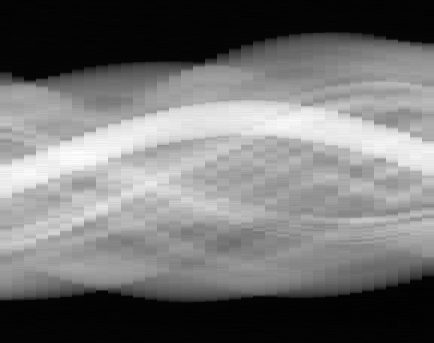}
         \caption{\footnotesize CT data}
         \label{ct data}
     \end{subfigure}
	\hspace{1em}
     \begin{subfigure}[t]{0.3\textwidth}
         \centering
         \includegraphics[width=\textwidth]{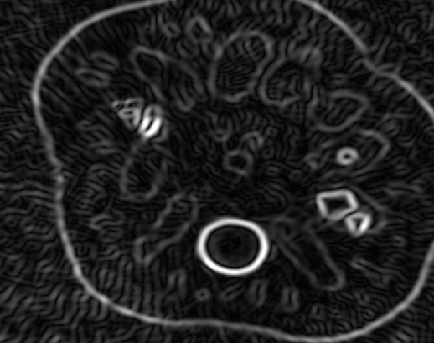}
         \caption{Method 1: $\abs{\nabla f_{\epsilon}}$, $\epsilon=3$}
         \label{method1:gmag}
     \end{subfigure}
   	\hspace{1em}
     \begin{subfigure}[t]{0.3\textwidth}
         \centering
         \includegraphics[width=\textwidth]{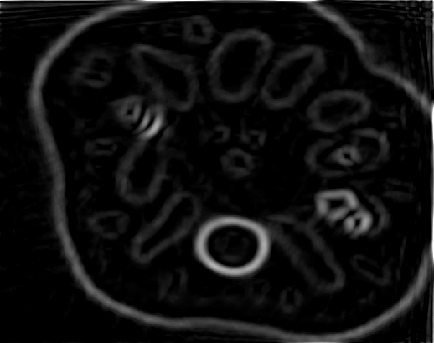}
         \caption{Method 2: $\abs{\nabla f_{\epsilon}}$, $\epsilon=6$, $\lambda=0.01$}
         \label{method2:gmag}
     \end{subfigure}
   
     \begin{subfigure}[t]{0.3\textwidth}
         \centering
         \includegraphics[width=\textwidth]{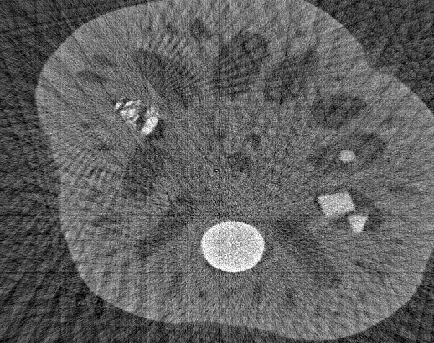}
         \caption{FBP reconstruction}
         \label{fbp rec}
     \end{subfigure}
     \hspace{1em}
     \begin{subfigure}[t]{0.3\textwidth}
         \centering
         \includegraphics[width=\textwidth]{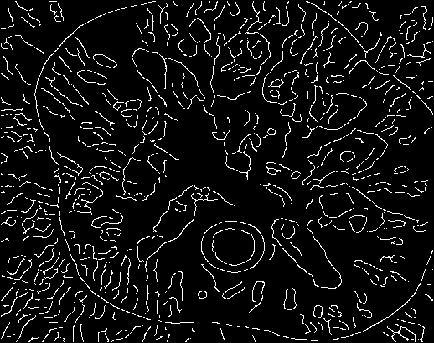}
         \caption{Method1: edge map (Canny)}
         \label{method1:E}
     \end{subfigure}
     \hspace{1em}
     \begin{subfigure}[t]{0.3\textwidth}
         \centering
         \includegraphics[width=\textwidth]{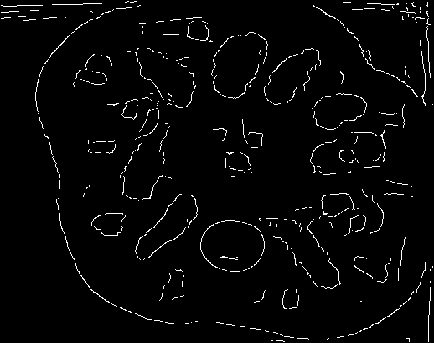}
         \caption{Method 2: edge map (Canny)}
         \label{method2:E}
     \end{subfigure}
\caption{Rebinned CT data of a lotus root \subref{ct data} and the corresponding FBP reconstruction \subref{fbp rec} from an angular range given by $[0,\pi)$ and $36$ evenly distributed angles, cf. \cite{FIPSlotus}. The magnitude of the gradient $\abs{\nabla f_{\epsilon}}$, in \subref{method1:gmag} and \subref{method2:gmag}, and the corresponding edge detection results using the Canny algorithm, see \subref{method1:E} and \subref{method2:E}.}
\label{fig:lotus edges}
\end{figure}

The CT data and the corresponding FBP reconstruction are shown in Fig. \ref{ct data} and \ref{fbp rec}. It can be clearly observed that the FBP reconstruction contains many undersampling artifacts (streaks) that could complicate the detection of edges in that image. Hence, we calculated the gradient directly from CT data using the methods 1 and 2, where we chose the parameters based on visual inspection of edge detection results. These are given in the caption of Fig. \ref{fig:lotus edges}.

In Fig. \ref{method1:gmag} and \ref{method1:E} one can see that the Gaussian smoothing cannot compensate for the undersampling artifacts and, thus, many edges in the edge map are not coming from actual image features. However, the $\ell^{1}$-regularization successfully reduces the number of artifacts and detects the edges more reliably, as can be seen in Fig. \ref{method2:gmag} and \ref{method2:E}. In other experiments we also observed that the method 2 outperforms the method 1 whenever the CT data was not sampled properly. For dense angular sampling, we found that both methods produce similar edge detection results.

\section{Discussion}

We presented two methods for calculating the gradient of a CT scan directly from CT data. As our first method, we introduced is a variant of the filtered backprojection algorithm and explained two different ways for its implementation. As our second approach, we introduced a sparsity based gradient recovery from CT data and showed in numerical experiments that this method is able to account for data undersampling and to provide more reliable edge detection results. Moreover, the second approach provides more flexibility and can be more easily applied to different scanning geometries.

\section{Acknowledgement}
The contribution by S.\, G. is part of a project that has received funding from the European Union’s Horizon 2020 research and innovation programme under the Marie Sk\l{}odowska-Curie grant agreement No 847476. The views and opinions expressed herein do not necessarily reflect those of the European Commission.

\bibliography{references.bib}

\begin{thebibliography}{1}

\bibitem{Natterer86}
Natterer F.
\newblock The mathematics of computerized tomography.
\newblock Stuttgart: B. G. Teubner; 1986.

\bibitem{FrikelQuinto2013}
Frikel J, Quinto ET.
\newblock Characterization and reduction of artifacts in limited angle
  tomography.
\newblock Inverse Problems. 2013;29(12):125007.

\bibitem{CannyOrig}
{Canny} J.
\newblock A Computational Approach to Edge Detection.
\newblock IEEE Transactions on Pattern Analysis and Machine Intelligence.
  1986;PAMI-8(6):679--698.

\bibitem{Lo2008}
Louis AK.
\newblock {Combining Image Reconstruction and Image Analysis with an
  Application to Two-Dimensional Tomography}.
\newblock SIAM J Img Sci. 2008;1:188--208.

\bibitem{Natterer01}
Natterer F, W{\"u}bbeling F.
\newblock Mathematical methods in image reconstruction.
\newblock SIAM Monographs on Mathematical Modeling and Computation.
  Philadelphia, PA: Society for Industrial and Applied Mathematics (SIAM);
  2001.

\bibitem{CScrt}
Candès EJ, Romberg JK, Tao T.
\newblock Stable signal recovery from incomplete and inaccurate measurements.
\newblock Comm Pure Appl Math. 2006;59(8):1207--1223.

\bibitem{FIPSlotus}
Bubba TA, Hauptmann A, Huotari S, et~al.
\newblock Tomographic X-ray data of a lotus root filled with attenuating
  objects. 2016; p. 1--11.
\newblock {arXiv:1609.07299 [physics{.}data-an]}.

\bibitem{DDD04}
Daubechies I, Defrise M, De~Mol C.
\newblock An iterative thresholding algorithm for linear inverse problems with
  a sparsity constraint.
\newblock Comm Pure Appl Math. 2004;57(11):1413--1457.

\end{thebibliography}

\end{document}